\documentclass{article}
\usepackage[]{epsf,epsfig,amsmath,amssymb,amsfonts,latexsym}
\usepackage[margin=1in]{geometry}
\usepackage{setspace}

\newtheorem{theorem}{Theorem}[section]
\newtheorem{lemma}[theorem]{Lemma}

\newtheorem{question}[theorem]{Question}

\newtheorem{proposition}[theorem]{Proposition}

\newtheorem{corollary}[theorem]{Corollary}
\newtheorem{definition}[theorem]{Definition}
\newtheorem{remark}[theorem]{Remark}

\def \PP {\mathbb P}
\def\qed{\hfill $\vcenter{\hrule height .3mm
\hbox {\vrule width .3mm height 2.1mm \kern 2mm \vrule width .3mm
height 2.1mm} \hrule height .3mm}$ \bigskip}
\def\P{\mathbb{P}}

\def\EE{\mathbb{E}}

\def\RR{\mathbb{R}}

\begin{document}
\title{Skorokhod Embeddings via Stochastic Flows on the Space of Gaussian Measures}
\author{Ronen Eldan}

\maketitle
\abstract{We present a new construction of a Skorohod embedding, namely, given a probability measure $\mu$ with zero expectation and finite variance, we construct an integrable stopping time $T$ adapted to a filtration $\mathcal{F}_t$, such that $W_T$ has the law $\mu$, where $W_t$ is a standard Wiener process adapted to the same filtration. We find several sufficient conditions for the stopping time $T$ to be bounded or to have a sub-exponential tail. In particular, our embedding seems rather natural for the case that $\mu$ is a log-concave measure and $T$ satisfies several tight bounds in that case. Our embedding admits the property that the stochastic measure-valued process $\{ \mu_t \}_{t \geq 0}$, where $\mu_t$ is as the law of $W_T$ conditioned on $\mathcal{F}_t$, is a Markov process. In view of this property, we will consider a more general family of Skorokhod embeddings which can be constructed via a kernel generating a stochastic flow on the space of measures. This family includes existing constructions such as the ones by Az\'ema-Yor (1979) and by Bass (1983), and thus suggests a new point of view on these constructions.}


\section{Introduction}
The \emph{Skorokhod embedding problem}, first presented by Skorokhod in \cite{S} (1961), was originally formulated as follows:
Given a prescribed centered probability measure $\mu$ whose second moment is finite and a standard Wiener process $W_t$ adapted to a filtration 
$\mathcal{F}_t$, can one find an integrable stopping time $T$, such that $W_T$ has the law $\mu$? \\

This problem has encouraged rather extensive research in the past 50 years (e.g., by Az\'ema, Bass, Dubins, Monroe, Ob\l\'{o}j, Root, Rost, Yor and many others), some of which is devoted to constructing new solutions, some to formulating and proving more general cases of this problem, and some to establishing certain properties of the existing solutions. A few examples of properties one would be interested to establish about a solution are bounds on moments of the stopping time $T$ given 
information about the measure $\mu$, monotonicity of $T$ with respect to quantities related to this measure and bounds related to the set $\{W_t, ~ 0 \leq t \leq T \}$. For an extensive review of many of these results, the reader is referred to \cite{Ob}. See also \cite{AHI} for a more recent construction of a solution. \\ 

The main contribution of the present note is to introduce a new solution to the Skorokhod embedding problem, based on the construction of a stochastic flow on the space of Gaussian measures. In addition, we will present several properties that this solution admits, mainly concerning bounds on the stopping time $T$ given additional assumptions on $\mu$. One of the advantages of this new solution is the existence of a formula with which the behaviour of $T$ can be analysed in many cases. In particular, this solution seems rather natural for log-concave measures, in which case we will derive sharp bounds for the moments of $T$. 

The construction has another property that may be notable: consider the measures $\mu_t$ defined to be the law of $W_T$ conditioned on $\mathcal{F}_t$. The process
$\{ \mu_t \}_{t \geq 0}$ is, in some sense, a Markov process whose "transition kernel" does not depend on the initial measure $\mu$. In view of this property, we will be able to relate our construction to other constructions in the literature and consider a more general framework for solutions to the Skorokhod problem.

Our construction is somewhat similar to the localization described in \cite{E} in the sense that the main mechanism behind it construction is a certain flow on the space of densities on $\RR$, defined by a system of stochastic differential equations. The central ideas behind it are described briefly in the beginning of Section 2. \\

Let us formulate our theorems. Throughout this note, $\mu$ will denote some fixed Borel probability measure on $\RR$. The only assumptions we will need for the construction of the Skorokhod embedding are that $\mu$ has expectation zero and a finite second moment:
\begin{equation} \label{unifint}
\int_{\RR} x^2 \mu(dx) < \infty ~ \mbox{  and  } ~ \int_{\RR} x \mu(dx) = 0.
\end{equation}
$$
~
$$
Let $W_t$ be a standard Wiener process adapted to a filtration $\mathcal{F}_t$. Our first goal will be to construct a stopping time $T_\mu$. The main properties of this stopping time are described in the following theorem, whose point is that $T_\mu$ induces a \emph{Skohorod embedding} of $\mu$ into the probability space of $\{W_t\}$. The actual definition of $T_\mu$ is postponed to the next section. 
\medskip
\begin{theorem} \label{mainthm}
Let $\mu$ be a measure on $\RR$ satisfying (\ref{unifint}). Then the stopping time $T_\mu$ satisfies the following properties: \\
(i) The event $\{ T_\mu \leq t\}$ is measurable with respect to the $\sigma$-algebra generated by $\{W_s \}_{s \leq t}$. \\
(ii) The random variable $W_{T_\mu}$ is distributed according to the law $\mu$. \\
(iii) One has $\EE[T_\mu] = Var[\mu]$.
\end{theorem}
\bigskip
The next theorems in this note establish bounds on the distribution of $T_\mu$ given that $\mu$ satisfies additional assumptions. Several estimates resembling some of our bounds have been established in \cite{AS} for a Skorokhod embedding based on the solution of a backwards stochastic differential equation. \\ 

A measure $\mu$ is said to be \emph{log-concave} if it is either a Dirac $\delta$-measure or has a density $f(x)$ with respect to the Lebesgue measure of the form $f(x) = e^{-\Phi(x)}$ where $\Phi:\RR \to \RR \cup \{ + \infty \}$ is a convex function. When the measure $\mu$ is log-concave, the stopping time $T_\mu$ admits a sub-exponential tail behaviour, namely we have the following.
\medskip
\begin{theorem} \label{thmlc}
There exist universal constants $c, C > 0$ such that, if $\mu$ is a log-concave measure with $\EE[\mu] = 0$ and $Var[\mu] = 1$, then
$$
\PP(T_\mu > t) < C e^{-ct}.
$$
\end{theorem}

\begin{remark}
The above result is tight up to the constants $c,C$. To see this, let $\mu$ be the measure whose density is $\frac{1}{2 \sqrt 2} e^{- |x| / \sqrt{2}}$,
and let $X$ have the law $\mu$. Then one has
$$
\PP(|X| > x) > 0.5 e^{-c_1 x}
$$
for some $c_1>0$. For each $x,t>0$ define the events
$$
A_t = \{T_\mu > t \}, ~~ B_{t,x} = \{\exists s < t \mbox{ such that } |W_s| > x \}.
$$
By combining the sub-Gaussian tail decay of the distribution of $W_t$ with Doob's martingale inequality (see \cite[Theorem II.1.7]{RY}), we have
$$
\PP(B_{t,x}) < C_2 e^{-c_2 \frac{x^2}{t}}
$$
for some $C_2,c_2>0$. Using a union bound, we have for all $t>0$,
$$
\PP(|X| > x) = \PP(|W_{T_\mu}| > x) \leq \PP(A_t) + \PP(B_{t,x}) 
$$
and therefore
$$
\PP(A_t) \geq 0.5 e^{-c_1 x} - C_2 e^{-c_2 x^2/t}.
$$
By optimizing over $x$ we get
$$
\PP(T_\mu > t) \geq C' e^{-c' t}, ~~ \forall t > 0
$$
for some constants $C',c'$.
\end{remark}
\bigskip
For log-concave measures whose modulus of log-concavity is bounded from below, the stopping time will be bounded according to the following tight estimate:
\begin{theorem} \label{thmunilc}
Let $\sigma > 0$. Let $\mu$ be a centered probability measure which admits a density satisfying
$$
\frac{\mu(dx)}{d x} = e^{- \frac{x^2}{2 \sigma^2} - \Phi(x)}
$$
where $\Phi:\RR \to \RR \cup \{ + \infty \}$ is a convex function. Then one has almost surely,
$$
T_\mu \leq \sigma^2.
$$
\end{theorem}
\begin{remark}
The above bound is tight, as demonstrated by the case that $\mu$ is a Gaussian measure whose variance is $\sigma^2$.
\end{remark}
\bigskip
For a Borel measure $\mu$ on $\RR$, we denote by $Supp(\mu)$ the support of $\mu$ which is the minimal closed set of full measure. Our next task is to address measures whose support is a compact set. For such measures we can give deterministic bounds on the stopping time if the measure is either log-concave or absolutely continuous with respect to the Lebesgue measure with density bounded between two constants in an interval. This is summarized in the next two theorems:
\begin{theorem} \label{thmcompartreg}
Let $\mu$ be a measure on $\RR$ which satisfies: \\
(i) $Supp(\mu)$ is an interval contained in $[-L,L]$ for some $L>0$. \\
(ii) $\mu$ is absolutely continuous with respect to the Lebesgue measure and $\alpha \leq \frac{\mu(dx)}{dx} \leq \beta$ for all $x \in Supp(\mu)$.
Then,
$$
T_\mu \leq 2 L^2 \frac{\beta}{\alpha}
$$
almost surely.
\end{theorem}
\bigskip
\begin{theorem} \label{thmcompartlc}
Let $\mu$ be a log concave measure on $\RR$ with $Supp(\mu) \subseteq [-L,L]$ for some $L>0$. Then,
$$
T_\mu \leq 2 L^2.
$$
almost surely.
\end{theorem}
\begin{remark}
The bound of the above theorem is tight up to the constant $2$. Indeed, let $\mu$ be the uniform measure over an interval of length $2L$. Then we have $\EE[T_\mu] = Var[\mu] = \frac{L^2}{3}$, therefore it cannot be the case that $T_\mu < \tfrac {L^2} 3$ almost surely.
\end{remark}
$$
~
$$
The structure of the remainder of this note is the following: In Section 2 we construct the stopping time $T_\mu$ and prove Theorem \ref{mainthm}. Section 3 deals with log-concave measures, in this section we prove theorems \ref{thmlc} and \ref{thmunilc}. In Section 4 we prove theorems \ref{thmcompartreg} and \ref{thmcompartlc}. In Section 5, we define the "Markov property" satisfied by our construction and relate our construction to existing constructions of Skorokhod embeddings in the context of this definition. In the appendix, we fill in some missing details left open in the construction of $T_\mu$. \\ 

Throughout this note, we use the following notation: for an It\^{o} process $X_t$, we denote by $d X_t$ the differential of $X_t$, and by $[X]_t$ the quadratic variation of $X_t$. For a pair of continuous time stochastic processes $X_t, Y_t$, the quadratic covariation will be denoted by $[X, Y]_t$. For a measure $\mu$ on $\RR$, we denote by $\EE[\mu]$ and $Var[\mu]$ its expectation and variance respectively. By $\frac{\mu(dx)}{dx}$ we denote the Radon-Nykom density of $\mu$ with respect to the Lebesgue measure at the point $x$ and likewise by $\frac{\mu(dx)}{\nu(dx)}$ we denote the density of $\mu$ with respect to $\nu$ at $x$. When we write $X \sim \mu$ we mean that the random variable $X$ is distributed according to the law of the measure $\mu$. As a convention, when stating that a relation between two random variables (e.g., equality or inequality) holds, we mean that this relation holds \emph{almost surely} unless stated otherwise.
\\ 

\emph{Acknowledgements}
I would like to thank Ofer Zeitouni and Joseph Lehec for fruitful discussions. Moreover, I am thankful to the two referees for reports containing extremely useful comments on a preliminary version of this note, and in particular for suggesting the connections between the definitions of Section 5 to the Bass and the Az\'ema-Yor embeddings, and to other known solutions.

\section{Construction of the embedding}

The goal of this section is to construct the stopping time $T_\mu$ and to establish some of its basic properties. \\ 

Let us briefly describe the idea behind our construction. Given a measure $\mu$, we will construct a random one-parameter family of measures, $\{\mu_t\}_{t \geq 0}$ such that $\mu_0=\mu$, and for which there exists some random time $T>0$ with $\EE[T] < \infty$ such that the density $\frac{\mu_t(dx)}{\mu(dx)}$ is a Gaussian density for all $0 \leq t < T$ and $\mu_t$ is a Dirac $\delta$-measure for $t \geq T$. Moreover, \\
(i) For any measurable $A \subset \RR$, the process $\mu_{t \wedge T}(A)$ is a martingale, \\
(ii) The process $W = \{ \int_{\RR} x \mu_t(dx) \}_{t \geq 0} $ is a Brownian motion for $0 \leq t < T$, \\
(iii) $\mu_t$ converges (in $L_2$) to a Dirac $\delta$-measure $\mu_T$ as $t \to T$, \\
(iv) The process $\{ \mu_t \}$ is an adapted process with respect to the filtration generated by the Brownian motion $W$. \\
By properties (i), (iii) and with the help of the optional stopping theorem we will deduce that for a measurable set $A$, one has $\EE[\mu_{T}(A)] = \mu(A)$. Recall that $T$ is the time in which $\mu_t$ becomes a Dirac measure which in turn tells us that $\P(\int_{\RR} x \mu_T(dx) \in A) = \EE[\mu_T(A)]$. In other words, $\int_{\RR} x \mu_T(dx)$ will have the law $\mu$. Property (ii) will ensure us that this quantity is, in fact, a Brownian motion taken at time $T$, and property (iv) will ensure us that no "extra randomness" is used. The construction of the densities $F_t(x) = \frac{\mu_t(dx)}{\mu(dx)}$ is best described in formula (\ref{contloc}) below, and the time $T=T_\mu$ will be defined as the time in which
the solution "explodes" (hence, the solution ceases to exist). 

The measure-valued process $\mu_t$ will be a essentially continuous time Markov process (see Section 5 for details). Loosely speaking, the transition rule of this Markov process will be such that given the measure $\mu_t$ and the Brownian increment $d W_t$, $\mu_{t+dt}$ will be the unique measure whose density with respect to $\mu_t$ is a linear function whose coefficients are chosen so that $\mu_{t+dt}$ is a probability measure and so that $\EE[\mu_{t+dt}] = \EE[\mu_t] + d W_t$. \\

We begin with some definitions. Let $\mu$ be a probability measure on $\RR$, satisfying (\ref{unifint}). For $c \in \RR$ and $b \geq 0$, we write
$$
V_{\mu}(b,c) = \int_{\RR} e^{c x - \frac{1}{2} b x^2 } \mu(dx)
$$
and define two functions,
$$
a_{\mu}(b,c) = V_{\mu}^{-1} (b,c) \int_{\RR} x e^{cx - \frac{1}{2} b x^2  }  \mu(dx),
$$
and
$$
A_{\mu}(b,c) = V_{\mu}^{-1} (b,c) \int_{\RR} (x - a_{\mu}(b,c))^2 e^{cx - \frac{1}{2} b x^2  }  \mu(dx).
$$
Note that $a_\mu(0,0) = \int_{\RR} x \mu(dx) = 0$ and $A_{\mu}(0,0) = \int_{\RR} x^2 \mu(dx)$. It is easy to verify that the assumption \eqref{unifint} implies that the functions $V_{\mu}$, $a_{\mu}$ and $A_{\mu}$ are smooth functions in the domain $(b,c) \in (0, \infty) \times \RR$. \\ 

Let $W_t$ be a standard Wiener process and consider the following system of stochastic differential equations:
\begin{equation} \label{stochastic1}
c_0 = 0, ~~ d c_t = A_{\mu}^{-1}(b_t, c_t) dW_t + A_{\mu}^{-2}(b_t, c_t) a_{\mu}(b_t, c_t) dt,
\end{equation}
$$
b_0 = 0, ~~ d b_t = A_{\mu}^{-2}(b_t, c_t) dt.
$$

First, we will explain why the solution exists under a stronger assumption, namely, that $\mu$ has some finite exponential moment. The proof for the more general case, assuming only that the second moment exists, is left for the appendix. Assume that there exists a constant $\delta > 0$ such that
$$
\int_{\RR} \left(e^{\delta x} + e^{- \delta x} \right ) \mu(dx) < \infty.
$$
Under this assumption, it is not hard to check that the functions $A_\mu^{-2} (b,c)$ and $a_\mu(b,c)$ are smooth functions on the set
$$
(b,c) \in \left ( (0,\infty) \times \RR \right ) \cup \left ( \{0\} \times [-\delta/2,\delta/2] \right ).
$$
Since $\frac{d}{dt} b_t |_{t=0} > 0$, these functions can be modified so that they are smooth on the set $[0, \infty) \times \RR$ without affecting 
the solution. In this case, we can use a standard existence and uniqueness theorem (see e.g., \cite{oksendal}, Section 5.2) to ensure the existence and uniqueness of a solution on some interval $[0,t_0)$ where $t_0$ is an almost-surely positive random variable. \\
\begin{remark}
Note that the fact that $\delta > 0$ is crucial for this argument, and the existence would not generally be true if $\delta = 0$. In general,
the functions $A_\mu, a_\mu$ may not be defined in any neighbourhood of $(0,0)$ of the form $[0, \epsilon] \times (-\epsilon, \epsilon)$, and are only bounded in parabolic sets the form $\{(b,c); \epsilon c^2 < b \}$. A-priori, in order to ensure the existence of the solution, one has to prove
that $(b_t,c_t)$ remain in such a set. This will be done in the appendix in an indirect fashion.
\end{remark}
$$
~
$$
We are now ready to define our stopping time $T_\mu$: it will be defined as the supremum over the set of times in which the solution to (\ref{stochastic1}) exists (and define $T_\mu=\infty$ if the solution exists for all $t \geq 0$). It is not hard to verify that the functions $A_\mu^{-2} (b,c)$ and $a_\mu(b,c)$ are smooth functions on any set in which $b$ is bounded away from zero. Since $b_t$ is increasing and since $\frac{d}{dt} b_t|_{t=0} > 0$, it follows that for any point $t$ in which $A_\mu^{-2} (b_t,c_t) > 0$ there exists some $\epsilon > 0$ such that the solution may be continued in the interval $[t, t+\epsilon]$. Consequently, we see that if one defines
\begin{equation} \label{defT}
T_\mu = \sup \{t \geq 0; ~ A_\mu(b_t,c_t) > 0 \},
\end{equation}
then almost surely, the solution of (\ref{stochastic1}) exists exactly in the interval $[0,T_\mu)$. Our next main goal is to show that $W_{T_\mu}$ has the law $\mu$. We abbreviate $T = T_\mu$.\\ 

We begin with the construction a 1-parameter family of measures $\mu_t$ by writing \\
\begin{equation} \label{defineF}
F_t(x) = V_{\mu}^{-1} (b_t, c_t) e^{c_t x - \frac{1}{2} b_t x^2  }
\end{equation}
and defining the measure $\mu_t$ by,
$$
\frac{\mu_t(dx)}{\mu(dx)}  = F_t(x)
$$
for all $0 \leq t < T$. Note that, by the above definitions, we have
$$
\mu_t(\RR) = \int_{\RR} F_t(x) \mu(dx) = V_{\mu}^{-1}  (b_t, c_t) \int_{\RR}  e^{c_t x - \frac{1}{2} b_t x^2  } \mu(dx) = 1
$$
hence $\mu_t$ is a probability measure. Also, abbreviate
$$
a_t = a_{\mu}(b_t, c_t), ~~ A_t = A_{\mu}(b_t, c_t), ~~ V_t = V_{\mu}(b_t, c_t).
$$
This gives by definition
$$
a_t = \int_{\RR} x \mu_t(dx), ~~ A_t = \int_{\RR} (x - a_t)^2 \mu_t(dx)
$$
hence $a_t$ and $A_t$ are respectively the expectation and the variance of the measure $\mu_t$. \\ 

The following lemma may shed some light on this construction.
\begin{lemma} \label{basicform}
For all $t \in [0, T)$ and for all $x \in \RR$, the process $F_t(x)$ satisfies the following set of equations:
\begin{equation} \label{contloc}
F_0(x) = 1, ~~ d F_t(x) = (x - a_t) A_t^{-1} F_t(x) d W_t,
\end{equation}
$$
a_t = \int_{\RR} x F_t(x) \mu(dx), ~~ A_t = \int_{\RR} (x - a_t)^2 F_t(x) \mu (dx).
$$
moreover, if $\xi(x)$ is a function satisfying $|\xi(x)| \leq C |x|^p$ for some constants $C,p>0$ then one has for all $t>0$,
\begin{equation} \label{stocfubini}
d \int_{\RR} \xi(x) \mu_t(dx) = \int_{\RR} \xi(x) d F_t(x) \mu(d x)
\end{equation}
\end{lemma}
\emph{Proof:} \\
Fix $x \in \RR$. We will show that $d F_t(x) = (x - a_t) A_t^{-1} F_t(x) d W_t$. The correctness of the other equations is obvious.
Define, 
$$
G_t(x) = V_t F_t(x) = e^{c_t x - \frac{1}{2} b_t x^2 }.
$$
Equation (\ref{stochastic1}) clearly implies that $[b]_t = 0$. Moreover,
$$
d c_t = A_t^{-1} dW_t + A_t^{-2} a_t dt
$$
It follows that,
$$
d [c]_t = A_t^{-2} dt.
$$
Using It\^{o}'s formula, we calculate
$$
d G_t(x) = \left ( x d c_t - \frac{1}{2} d b_t x^2 + \frac{1}{2} x^2 d [c]_t \right ) G_t(x)
$$
$$
= \left ( x A_t^{-1} dW_t + x A_t^{-2} a_t dt - \frac{1}{2} A_t^{-2} x^2 dt + \frac{1}{2} A_t^{-2} x^2 dt \right ) G_t(x)
$$
$$
= x \left ( A_t^{-1} dW_t + A_t^{-2} a_t dt \right ) G_t(x).
$$
Next, we claim that we may apply a stochastic Fubini theorem to get
$$
d \int_{\RR} G_t(x) \mu(dx) = \int_{\RR} d G_t(x) \mu(dx).
$$
Indeed, since $G_t(x)^2 \leq 1$ almost surely for all $x,t$, it is easy to verify that the conditions of \cite[Theorem 2.2]{V} hold and the last equation follows. Therefore, we can calculate
$$
d V_t = d \int_{\RR} e^{c_t x - \frac{1}{2} b_t x^2} \mu(dx)
$$
$$
=\int_{\RR} d G_t(x) \mu(dx) = \int_{\RR} x \left ( A_t^{-1} dW_t + A_t^{-2} a_t dt \right ) G_t(x) \mu(dx)
$$
$$
= V_t a_t \left ( A_t^{-1} dW_t + A_t^{-2} a_t dt \right ).
$$
So, using It\^{o}'s formula again,
$$
d V_t^{-1} = - \frac{d V_t}{V_t^2} + \frac{d [V]_t}{V_t^3}
$$ 
$$
= - V_t^{-1} a_t \left ( A_t^{-1} dW_t + A_t^{-2} a_t dt \right ) + V_t^{-1} A_t^{-2} a_t^2 dt.
$$
Applying It\^{o}'s formula one last time yields,
$$
d F_t(x) = d (V_t^{-1} G_t(x)) 
$$
$$
= G_t(x) d V_t^{-1} + V_t^{-1} d G_t(x) + d [V^{-1}, G(x)]_t
$$
$$
= - V_t^{-1} a_t \left (  A_t^{-1} dW_t + A_t^{-2} a_t dt \right ) G_t(x) + V_t^{-1} A_t^{-2} a_t^2  G_t(x) dt + 
$$
$$
+ V_t^{-1} x ( A_t^{-1} dW_t + A_t^{-2} a_t dt ) G_t(x) - A_t^{-2} a_t x V_t^{-1} G_t(x) dt
$$
$$
= A_t^{-1} dW_t ( x - a_t ) F_t(x).
$$
This finishes the proof of formula \eqref{contloc}. In order to prove \eqref{stocfubini}, we recall that that for any $t>0$ one has $b_t > 0$. Consider the expression
$$
\psi_t(x) := \xi(x) A_t^{-1} ( x - a_t ) F_t(x) 
$$
and fix some $t_0 > 0$. Clearly, since $|\xi(x)| < C |x|^p$ and since $F_{t_0}$ has a sub-Gaussian tail, we have that $\psi_{t_0}(\cdot)$ is bounded by some constant, say $M$. Define the stopping time
$$
\tau = \inf \{s>t_0; ~ \exists x_0 \mbox{ such that } \psi_t(x_0) \geq 10M \}.
$$
By definition, $\psi_t(x)$ is bounded in the interval $t_0 < t < \tau$, and we can invoke a standard stochastic Fubini theorem (e.g., \cite[Theorem 2.2]{V}) to get
$$
\int_{t_0}^{\tau \wedge s} \int_{\RR} \psi_t(x) \mu(dx) d W_t = \int_{\RR} \int_{t_0}^{\tau \wedge s} \psi_t(x) d W_t \mu(dx).
$$
By continuity we have $\tau > t_0$ almost surely, so by taking the limit $s \to t_0$ we conclude that for all $t>0$ one has
$$
d \int_{\RR} \xi(x) F_t \mu(dx)  = \int_{\RR} \xi(x) d F_t(x) \mu(dx)
$$
and equation \eqref{stocfubini} follows. The lemma is complete.
\qed \\

At this point, we would like to formally extend the definition of $\mu_t$ to $t \geq T$. To that end, fix a measurable set $B \subset \RR$ and consider the 
process $\mu_t(B)$. According to the above lemma we have
$$
d \mu_t(B) = \int_B d F_t(x) \mu(dx) = A_t^{-1} \int_B (x - a_t) F_t(x) \mu(dx) d W_t
$$
hence, this process is a local martingale. Recall that $\mu_t$ is almost surely a probability measure, which implies that $\mu_t(B) \in [0,1]$. The martingale convergence theorem thus implies that the limit
$$
\mu_T(B) := \lim_{t \to T^-} \mu_t(B)
$$
almost surely exists. Finally, we define for all $t>T$, $\mu_t := \mu_T$. Note that, at this point, $\mu_T$ is a set function defined on measurable sets (later on we will establish that it is a measure).  Now, since every bounded local martingale is also a martingale, we establish the following result
\begin{lemma} \label{lemmartingale}
Let $B \subset \RR^n$ be measurable. Then $\mu_{t \wedge T} (B)$ is a martingale.
\end{lemma}
\medskip
The next lemma is a simple calculation that extracts one of the main points of the construction: the process of the centers of mass of $\mu_t$ is a Brownian motion.
\begin{lemma} \label{lemdat}
For all $0 \leq t < T$, one has $a_t = W_t$.
\end{lemma}
\emph{Proof:} \\
Using formulas (\ref{contloc}) and \eqref{stocfubini}, we calculate
$$
d a_t = d \int_{\RR} x F_t(x) \mu(dx) = \int_{\RR} x d F_t(x) \mu(dx)
$$ 
$$
= \left (\int_{\RR} x (x - a_t) A_t^{-1} F_t(x) \mu(dx) \right ) dW_t.
$$
Now, by the definition of $a_t$, one has
$$
\int_{\RR} a_t (x - a_t) A_t^{-1} F_t(x) \mu(dx) = A_t^{-1} \left (a_t \mu_t(\RR) - \int_{\RR} x \mu_t(dx) \right ) = 0
$$
Joining the two previous equations together gives
\begin{equation} \label{dat}
d a_t = \left (\int_{\RR} (x - a_t)^2 A_t^{-1} F_t(x) \mu(dx) \right ) dW_t = 
\end{equation}
$$
A_t^{-1} \left (\int_{\RR} (x - a_t)^2 F_t(x) \mu(dx) \right ) dW_t = d W_t.
$$
\qed \\ 

Next, we prove
\begin{lemma} \label{Tfinite}
One has for all $0 < t < T$,
\begin{equation} \label{dAt2}
d A_t = 
\left (\int_{\RR} (x - W_t)^3 \mu_t(dx) \right ) A_t^{-1} d W_t - dt.
\end{equation}
Moreover, $T$ is almost-surely finite.
\end{lemma}
\emph{Proof:} \\
Using formulas \eqref{contloc} and \eqref{stocfubini}, we have for all $t>0$,
$$
d A_t = d \left ( \int_{\RR} x^2 F_t(x) \mu(dx) - a_t^2 \right ) = d \left ( \int_{\RR} x^2 F_t(x) \mu(dx) - W_t^2 \right )
$$
$$
= \left ( \int_{\RR} x^2 A_t^{-1} (x - W_t) \mu_t(dx) \right ) d W_t - 2 W_t d W_t - dt
$$
(By definition of $A_t$ and since $\int W_t (x-W_t) \mu_t(dx) = 0$)
$$
= \left ( \int_{\RR} x^2 (x - W_t) \mu_t(dx) \right )  A_t^{-1} d W_t - 2 W_t \left (\int_{\RR} x (x-W_t) \mu_t(dx) \right )  A_t^{-1} d W_t - dt
$$
$$
= \left ( \int_{\RR} (x^2 - 2 W_t x + W_t^2) A_t^{-1} (x - W_t) \mu_t(dx) \right ) d W_t - dt
$$
$$
= \left (\int_{\RR} (x - W_t)^3 \mu_t(dx) \right ) A_t^{-1} d W_t - dt
$$
which settles (\ref{dAt2}). \\ 

To see that $T$ is almost surely finite, write $X_t = A_t + t$ and $X_T = T$. The above equation, along with the fact that $\lim_{t \to T^-} A_t = 0$, suggests that $X_{t \wedge T}$ is a local-martingale, and since it is bounded from below, it is also a supermartingale. Suppose by contradiction that with positive probability, a solution exists for all $t>0$. This implies that $A_t$ exists and is positive for all $t>0$. By the martingale convergence theorem, we have 
$$\PP \left . \left (\lim_{t \to \infty} X_t \mbox{ exists} \right | ~ \forall t>0, X_t \geq 0 \right) = 1,$$
but observe that when $\lim_{t \to \infty} X_t$ exists then $\lim_{t \to \infty} A_t = - \infty$ which is clearly impossible. The lemma is complete.
\qed \\ \\
We are finally ready to prove that $T_\mu$ is a Skorokhod embedding. \\ \\
\emph{Proof of Theorem \ref{mainthm}:} \\
Part (i) of the theorem is obvious from the definition of $T$. To prove part (ii), let $\varphi(x)$ be a smooth, compactly supported function. We have for all $0 \leq t < T$,
$$
\left | \int_{\RR} \varphi(x) \mu_t(dx) - \varphi(W_t) \right | = \left | \int_{\RR}  (\varphi(x) - \varphi(a_t)) \mu_t(dx) \right | \leq
$$
$$
\sup_{x \in \RR} |\varphi'(x)| \int_{\RR} |x - a_t| \mu_t(dx) \leq \sup_{x \in \RR} |\varphi'(x)| A_t^{1/2},
$$
where in the first passage we used Lemma \ref{lemdat} and in the last passage we used the Cauchy-Schwartz inequality and the fact that $Var[\mu_t] = A_t$. Since $\lim_{t \to T^-} A_t = 0$, it follows that
\begin{equation} \label{limtoT}
\lim_{t \to T^-} \int \varphi(x) \mu_t(dx) = \varphi(W_T),
\end{equation}
Therefore $\mu_T$ is a Dirac probability measure with an atom at $W_T$ (and as promised above, we have established that it is a measure). Recall that for a measurable set $B$, $\mu_{t \cap T}(B)$ is a martingale (Lemma \ref{lemmartingale}). It follows that $\int \varphi(x) \mu_t(dx)$ is also a martingale. By the optional stopping theorem,
\begin{equation} \label{stoppingeq}
\EE \left . \left [ \int_{\RR} \varphi(x) \mu_{T \wedge t} (dx) \right | ~ \mathcal{F}_{s \wedge T} \right ] = \int_{\RR} \varphi(x) \mu_{T \wedge s}(dx), ~~ \forall t > s \geq 0.
\end{equation}
Now, since $\varphi$ is bounded and since $\mu_{t \wedge T}$ is a probability measure, one has
$$ 
\int_{\RR} |\varphi(x)| \mu_{T \wedge t} (dx) \leq \sup_{x \in \RR} |\varphi(x)|, ~~ \forall t \geq 0
$$
and the dominated convergence theorem implies that
$$
\lim_{t \to \infty} \EE \left . \left [ \int_{\RR} \varphi(x) \mu_{T \wedge t} (dx) \right | ~ \mathcal{F}_{s \wedge T} \right ] = \EE \left . \left [ \lim_{t \to T^-} \int_{\RR}  \varphi(x) \mu_{t} (dx) \right | ~ \mathcal{F}_{s \wedge T} \right ] 
$$
$$
= \EE \left . \left [\int_{\RR}  \varphi(x) u_{T} (dx) \right | ~ \mathcal{F}_{s \wedge T} \right ]
$$
for all $s \geq 0$ where we used the fact that $T$ is almost surely finite, proven in Lemma \ref{Tfinite}, with formula (\ref{limtoT}). Combining the last equality with (\ref{stoppingeq}) gives,
\begin{equation} \label{stoppingeq2}
\EE \left . \left[\int_{\RR} \varphi(x) u_T(dx) \right | ~ \mathcal{F}_{s \wedge T}  \right ] = \int_{\RR} \varphi(x) \mu_{s \wedge T}(dx)
\end{equation}
for all $s \geq 0$. Taking $s=0$ proves part (ii) of the theorem. \\ 

To prove part (iii) of the theorem, observe that it follows from the optional stopping theorem that
\begin{equation} \label{stoppingsq}
\EE[ W_{T \wedge t}^2 ] = \EE[T \wedge t], ~~ \forall t > 0.
\end{equation}
By taking the limit $t \to \infty$ on both sides, we see that it suffices to show that
\begin{equation} \label{stoppingsq2}
\lim_{t \to \infty} \EE[ W_{T \wedge t}^2 ] = \EE[W_T^2] = Var[\mu].
\end{equation}
To that end, for all $t>0$ define $X_t = W_T - W_{t \wedge T}$. Equation (\ref{stoppingeq2}) implies that for any compactly supported continuous test function $\varphi$, one has
$$
\EE[\varphi(X_s) | \mathcal{F}_{s \wedge T}] = \int_{\RR} \varphi(x - W_{s \wedge T}) u_{s \wedge T}(dx), ~~ \forall s \geq 0,
$$
where we used the fact that $W_{s \wedge T}$ is $\mathcal{F}_s$-measurable. By considering a monotone increasing
sequence $\varphi_n(x)$ of positive, compactly supported functions satisfying $\lim_{n \to \infty} \varphi_n(x)  = x^2$ for all $x \in \RR$, the monotone convergence theorem together with the last equation yield
$$
Var[X_s | \mathcal{F}_{s \wedge T}] = \lim_{n \to \infty} \EE[\varphi_n (X_s) | \mathcal{F}_{s \wedge T}] 
$$
$$
= \lim_{n \to \infty} \int_{\RR} \varphi_n (x - W_{s \wedge T}) u_{s \wedge T}(dx) =  A_{s \wedge T}, ~~ \forall s \geq 0
$$
(where we define $A_T = 0$). Moreover, since $\int_{\RR} x u_{s \wedge T}(dx) = W_{s \wedge T}$, it follows that $Cov(X_s, W_{s \wedge T}) = 0$ for all $s \geq 0$, which gives
$$
Var[X_s] + Var[W_{s \wedge T}] = Var[\mu]
$$
for all $s \geq 0$. Consequently,
$$
Var[W_{t \wedge T}] \leq Var[\mu], ~~ \forall t \geq 0
$$
which implies that
$$
\lim_{t \to \infty} \EE[ W_{T \wedge t}^2 ] \leq Var[\mu]
$$
where the fact that the limit exists follows from the fact that the limit of the right hand side of formula (\ref{stoppingsq}) exists. We conclude
that
$$
\EE[T] = \lim_{t \to \infty} \EE[T \wedge t] \leq Var[\mu].
$$
Now, since $T$ has a finite first moment, we can use the optional stopping theorem once again with the martingale $W_t^2 - t$
to get
\begin{equation} 
\EE[ W_{T}^2 ] = \EE[T]
\end{equation}
and the theorem is complete.
\qed

\section{Log concave measures}
We begin by recalling a few basic things about log-concave measures. A log concave measure $\mu$ on $\RR$ is either a Dirac measure or is absolutely continuous with respect to the Lebesgue measure. 

A central tool we will use will be the following well-known estimate, proven via integration by parts:
\begin{theorem}
Let $V: \RR \to \RR$ be a strictly convex function, such that $\int_{\RR } e^{-V(x)} dx = 1$. Let $\mu$ be a probability measure on 
$\RR$ defined by $\frac{\mu(dx)}{dx} = e^{- V(x)}$. Then for every smooth function $f:\RR \to \RR$,
$$
\int_{\RR} \left (f (x)- \int_{\RR} f(x) \mu(dx) \right )^2 \mu(dx) \leq \int_{\RR} \left ( V'' (x) \right )^{-1} (f' (x))^2  \mu (dx).
$$
\end{theorem}
\medskip
\begin{remark}
The above theorem is merely the one-dimensional version of a theorem of Brascamp-Lieb from \cite{BL}.
\end{remark}
\medskip
An application of this theorem with the function $f(x) = x$ gives,
\begin{proposition} \label{bakryemery}
Let $\phi: \RR \to \RR \cup \{\infty\}$ be a convex function and let $\sigma>0$.
Suppose that $\mu$ is a measure satisfying
$$
\mu(dx) = Z e^{-\phi(x) - \frac{1}{2 \sigma^2} |x|^2} dx
$$
with $Z>0$ being a normalizing constant so that $\mu$ is a probability measure. Then one has,
$$
Var[\mu] \leq \sigma^2.
$$
\end{proposition}
$$
~
$$
As a corollary, we have the following:
\begin{corollary} \label{corrlc}
If $\mu$ is a log-concave measure then
\begin{equation} \label{Atbound}
A_t \leq b_t^{-1}
\end{equation}
for all $0 \leq t < T$.
\end{corollary}
\emph{Proof:}
Recall the formula (\ref{defineF}), and apply Proposition \ref{bakryemery} with $\sigma^2 = \frac{1}{b_t}$.
\qed
$$
~
$$
Recall that $db_t = A_t^{-2} dt$. In light of this equation, and with the help of the above corollary, we have the following bound for $T$:

\begin{lemma} \label{explboundlem}
There exists a universal constant $C>0$ such that the following holds whenever $\mu$ is a log-concave measure: define the stopping time 
$$
\tau = \min(\inf \{t; ~ A_t \geq 2 \}, 1).
$$ 
One has almost surely,
$$
T \leq C + \frac{C}{\tau}.
$$
\end{lemma}
\emph{Proof:}
If $T \leq 1$ then we're done. Otherwise, recall that 
\begin{equation} \label{defbtrecall}
d b_t = A_t^{-2} dt
\end{equation}
we note that by the definition of $\tau$ and by this equation,
$$
b_\tau = \int_0^\tau d b_t = \int_0^\tau A_t^{-2} dt \geq \int_0^\tau \frac{1}{4} dt \geq \frac{\tau}{4}.
$$
Combine \eqref{defbtrecall} with equation (\ref{Atbound}) to get,
$$
\frac{d}{dt} b_t = A_t^{-2} \geq b_t^2.
$$
Consequently, for all $\tau < t < T$,
$$
\frac{1}{b_t} = \frac{1}{b_\tau} - \int_\tau^t \frac{\frac{d}{dt} b_t}{b_t^2} dt \leq \frac{1}{b_\tau} - (t - \tau) \leq \frac{4}{\tau} - (t - \tau)
$$
and since $\frac{1}{b_t} \geq 0$, we conclude that
$$
T \leq \tau + \frac{4}{\tau} \leq 1 + \frac{4}{\tau}.
$$
\qed \\ \\
We are now ready to prove that when $\mu$ is log-concave, $T_\mu$ has a sub-exponential tail. \\
\emph{Proof of Theorem \ref{thmlc}:} \\
Recall equation (\ref{dAt2}),
\begin{equation}
d A_t = 
\left (\int_{\RR} (x - W_t)^3 \mu_t(dx) \right ) A_t^{-1} d W_t - dt
\end{equation}
Define $S_t = \left (\int_{\RR} (x - W_t)^3 \mu_t(dx) \right ) A_t^{-1}$. A well-known fact about isotropic log-concave measures (see for example \cite[Lemma 5.7]{LV}) is that for every
$p \geq 2$ there exists a constant $c(p)$ such that for every log-concave measure $\nu$ on $\RR$,
$$
\int_{\RR} |x - \EE[\nu] |^p \nu(dx) \leq c(p) Var[\nu]^{p/2}.
$$
Using the above with the measure $\mu_t$ and with $p=3$ gives $| S_t | \leq C_1 A_t^{1/2}$ for some universal constant $C_1 > 0$ and for all $t>0$. With the definition of $\tau$, this gives 
\begin{equation} \label{Stbound}
S_t < 2 C_1, ~~ \forall 0 \leq t < \min(\tau,T).
\end{equation}
Next, define
$Y_t = A_t + t - 1$. By (\ref{dAt2}), we learn that $Y_t$ is a martingale. By the Dambis / Dubins-Schwartz theorem, there exists a monotone time change $\Theta(t)$ such that $Y_{\Theta(t)} \sim \tilde W_t$ where $\tilde W_t$  is a standard Wiener process defined in the interval $[0, \Theta^{-1}(T))$. Moreover, 
$$
\Theta'(t) = \left ( \frac{d}{dt} [Y]_t \right)^{-1} = S_t^{-2}.
$$
Equation (\ref{Stbound}) implies, 
\begin{equation}
\Theta(t) \geq c_2 t, ~~ \forall 0 \leq t \leq \min(\tau,T).
\end{equation}
for some universal constant $c_2 > 0$. An application of the so-called reflection principle now gives,
$$
\PP \left (\max_{t \in [0, p] } \tilde W_t \geq 2 \right )
$$
$$
= 2 \PP (\tilde W_p \geq 2 ) < C_3 e^{- 1 / p }.
$$
for some universal constant $C_3 > 0$, which implies that
$$
\PP \left (\frac{1}{\tau} > s \right ) \leq \PP \left (\max_{0 \leq t \leq s^{-1}} Y_t > 2 \right ) \leq
$$
$$
\PP \left (\max_{0 \leq c_2 t \leq s^{-1}} Y_{\Theta(t)} > 2 \right ) < C_3 e^{- c_2 s}.
$$
combining the last equation with Lemma \ref{explboundlem} finishes the proof.
\qed \\ 
We move on to proving Theorem \ref{thmunilc}, which states that if the density of $\mu$ with respect to some Gaussian measure is a log-concave
function, then $T$ is bounded by the variance of this Gaussian measure. \\ \\
\emph{Proof of Theorem \ref{thmunilc}:} \\
Thanks to the assumption of the theorem and to equation (\ref{defineF}), we know that for all $0 \leq t < T$, $\mu_t$ has the form
$$
\frac{\mu_t(dx)}{dx}  = e^{-\left ( \frac{b_t}{2} + \frac{1}{2 \sigma^2} \right ) x^2 - \Phi_t(x) }
$$
for some $\Phi_t:\RR \to \RR \cup \{\infty \}$ convex. Along with Proposition \ref{bakryemery}, this gives
$$
A_t^{-1} \geq \sigma^{-2} + b_t, ~~ \forall 0 \leq t < T.
$$
Define $e_t = \sigma^{-2} + b_t$. Combine this with the equation defining $b_t$, $d b_t = A_t^{-2} dt$ and with equation (\ref{Atbound}) to get,
$$
\frac{d}{dt} e_t = A_t^{-2} \geq e_t^2.
$$
Therefore, for all $t < T$ one has
$$
\frac{1}{e_t} = \frac{1}{e_0} - \int_0^t \frac{ \frac{d}{dt} e_t  }{ e_t^2} dt \leq  \frac{1}{\sigma^2} - t.
$$
Since $\frac{1}{e_t} > 0$ for all $t < T$, this implies 
$$
T \leq \sigma^2.
$$
\qed

\section{Measures with bounded support}
Let $\mu$ be a measure supported in the interval $[-L,L]$. Our main mean of using this fact will be the obvious observation that
\begin{equation} \label{obvious}
A_t \leq L^2, ~~ \forall 0 \leq t < T.
\end{equation}
The next lemma will be the main ingredient allowing us to take advantage of the fact that a measure has a density bounded between two constants on its support: \\
\begin{lemma} \label{AtBdDen}
Suppose that $\mu$ is absolutely continuous with respect to the Lebesgue measure on $\RR$ supported on some interval $I$ and suppose that $f(x) = \frac{\mu(dx)}{dx}$ satisfies 
$$
0 < \alpha \leq f(x) \leq \beta, ~~ \forall x \in I
$$
Let $a \in \RR$ and $b > 0$ and let $\nu$ be a probability measure defined by the equation
$$
\frac{\nu(dx)}{\mu(dx)} = Z^{-1} e^{- \frac{b}{2} (x - a) ^2}
$$
where $Z>0$ is the normilizing constant. Then
$$
Var[\nu] \leq \frac{\beta}{\alpha b}.
$$
\end{lemma}
\emph{Proof:} \\
Define
$$
x_0 = \frac{\int_{I} x e^{- \frac{b}{2} (x - a) ^2} dx }{\int_{I} e^{- \frac{b}{2} (x - a) ^2} dx }.
$$
An application of Proposition (\ref{bakryemery}) with the function $\mathbf{1}_{ \{x \in I \} } e^{- \frac{b}{2} (x - a) ^2}$ gives
$$
\frac{\int_{I} (x - x_0)^2 e^{- \frac{b}{2} (x - a) ^2} dx  }{\int_{I} e^{- \frac{b}{2} (x - a) ^2} dx } \leq \frac{1}{b}.
$$
Now, we have
$$
Var[\nu] \leq \int_{I} (x - x_0)^2 \nu(dx) = \frac{ \int_{I} (x - x_0)^2 e^{- \frac{b}{2} (x - a) ^2} \mu(dx)}{ \int_{I} e^{- \frac{b}{2} (x - a) ^2} \mu(dx)} \leq 
$$
$$
\frac{\beta}{\alpha} \frac{ \int_{I} (x - x_0)^2 e^{- \frac{b}{2} (x - a) ^2} dx}{ \int_{I} e^{- \frac{b}{2} (x - a) ^2} dx} \leq \frac{\beta}{\alpha b}.
$$
\qed \\ 
 
By combining the above lemma with equation (\ref{obvious}), we establish the bound for measures supported on an interval whose density
is bounded between two constants: \\ \\
\emph{Proof of Theorem \ref{thmcompartreg}:} \\
We conclude from the previous lemma that for all $0 \leq t < T$,
\begin{equation} \label{Atsmall22}
A_t \leq \frac{\beta}{\alpha} b_t^{-1}.
\end{equation}
Using this estimate and the estimate (\ref{obvious}) with the equation (\ref{stochastic1}) for $db_t$ 
\begin{equation} \label{eq1112}
\frac{d}{dt} b_t = A_t^{-2} \geq \max \left ( \frac{\alpha^2}{\beta^2} b_t^2, \frac{1}{L^4} \right ).
\end{equation}
Let $g(x)$ be a function satisfying
$$
g'(0) = \frac{1}{L^4}, ~~ g'(x) = \frac{\alpha^2}{\beta^2} g(x)^2
$$
Then,
$$
g(x) = \frac{1}{\frac{L^2 \alpha}{\beta} - \frac{\alpha^2}{\beta^2} x }.
$$
Define $t_0 = \frac{L^2 \beta}{\alpha}$ and note that, by (\ref{eq1112}),
$$
b_{t_0} \geq \frac{\beta}{\alpha L^2} = g(0).
$$
By a standard comparison theorem
$$
b_{t + t_0} \geq g(t), ~~ \forall 0 \leq t < T,
$$
which implies that for some $t_1 \leq 2 L^2 \frac{\beta}{\alpha}$, one has
$$
\lim_{t \to t_1} b_t = + \infty
$$
In light of formula (\ref{Atsmall22}), this implies that
$$
T \leq 2 L^2 \frac{\beta}{\alpha}
$$
and the proof is complete.
\qed \\ 

The proof of Theorem \ref{thmcompartlc} follows the same lines, only Lemma \ref{AtBdDen} is replaced by Proposition \ref{bakryemery}: \\ \\
\emph{Proof of Theorem \ref{thmcompartlc}:} \\
Using Corollary \ref{corrlc} and equation (\ref{obvious}) gives,
$$
A_t^{-2} \geq \max \left (b_t^2, \frac{1}{L^4} \right ).
$$
Plugging this into the formula for $db_t$, equation (\ref{stochastic1}), gives
$$
\frac{d}{dt} b_t \geq \max \left (b_t^2, \frac{1}{L^4} \right ).
$$
now follow the proof of Theorem \ref{thmcompartreg}, noting that equation (\ref{eq1112}) holds with $\frac{\alpha}{\beta} = 1$.
\qed

\section{Embeddings with a Markov property}
The goal of this section is to suggest a new point of view on some existing solutions to the Skorokhod embedding problem, including our solution. This point of view is related to another notable property admitted by our construction, namely, that the process defined by considering the measure $\mu$ conditioned on the filtration $\mathcal{F}_t$ is Markovian. By introducing two definitions related to this property, we will be able to view our construction as a member of a more general family of solutions to the Skorokhod problem; this is a family of solutions which satisfy a stochastic flow equation similar to \eqref{contloc}, each member of which is associated with a different kernel. We will show that this family contains two known solutions from the literature, constructed by Az\'ema-Yor \cite{AY} and by Bass \cite{Bass}. \\

Let $W_t$ be a standard Wiener process with a corresponding filtration $\mathcal{F}_t$. Let $\mathcal{M}$ be the space of Borel probability measures $\mu$ on $\RR$ such that $Var[\mu] < \infty$ and let $\mathcal{M}' \subset \mathcal{M}$ be the subset of measures whose centroid is the origin. Let $\mathcal{T}$ be the space of $\mathcal{F}$-stopping times. We define a \emph{Skorokhod embedding scheme} as a function $T:\mathcal{M}' \to \mathcal{T}$ taking $\mu$ to a stopping time $T(\mu)$ such that $W_{T(\mu)}$ has law $\mu$. In the following, we will abbreviate $T=T(\mu)$ when clarity is not affected.  \\

For a measure $\mu$, a Skorokhod embedding scheme $T$ and any $t>0$, define a random measure $\mu_t = \mu_t[\mu,T]$ by
\begin{equation} \label{defmutt}
\mu_t(E) = \mu_t[\mu,T](E) := \PP(W_{T(\mu)} \in E | ~ \mathcal{F}_t )
\end{equation}
for all measurable $E \subset \RR$. It is easy to verify that $T(\mu)$ is uniquely determined by the measure $\mu_t$ using the formula
\begin{equation} \label{Tmut}
T(\mu) = \sup \{t|~ Var[\mu_t] > 0\}.
\end{equation}
$$ ~ $$
For a measure $\nu \in \mathcal{M}$ and $s \in \RR$, denote by $L_{s} \nu$ the measure $\nu$ translated by $s$. Also define $\mathrm{C}(\nu) := L_{-\EE[\nu]} \nu \in \mathcal{M}'$, hence the measure $\nu$ translated by $-\EE[\nu]$ so that it becomes centered. Our first definition reads
\begin{definition} (Markov property).
We will say that a Skorokhod embedding scheme $T$ has the \emph{Markov} property if for every measure $\mu \in \mathcal{M}'$ and every $t>0$, we have almost surely that
$$
\left . \bigl (\mu_{t+s}[\mu,T] \bigr )_{s \geq 0} \right | \mathcal{F}_t \sim \bigl (L_{\EE[\nu]} \mu_s[\mathrm{C}(\nu),T] \bigr)_{s \geq 0}
$$
where in the right hand side, the measure $\nu$ is chosen to be the measure $\mu_t[\mu,T]$. In other words, almost surely the sequence of measures $\mu_{t+s}[\mu,T]$ (as a sequence parametrized by $s$), conditioned of $\mathcal{F}_t$, has the same distribution as the sequence of measures $\mu_s [\mathrm{C}(\nu),T]$, with the initial measure $\nu$ being $\mu_t$, up to translations.
\end{definition}

Roughly speaking, an embedding scheme has the Markov property if $\mu_{t + dt}$ depends only on $\mu_t$ and on $d W_t$. Next, we denote by $\mathcal{MF}$ the space of measurable functions $f:\RR \to \RR$. The main definition we are interested in is:
\begin{definition} (Analytic Markov property)
A Skorokhod embedding scheme $T$ has the analytic \emph{Markov} property if there exists a function $K:\mathcal{M} \to \mathcal{MF}$ such that for every $\mu \in \mathcal{M}$ the associated stopping time $T$ and the family of probability measures $\{\mu_t\}_{t>0}$ defined in \eqref{defmutt} satisfy
\begin{equation} \label{markoveq}
\mu_t(A) = \mu_0(A) + \int_0^t \left (\int_A K[\mu_s](x) \mu_s(dx) \right ) d W_s.
\end{equation}
for all measurable $A \subset \RR$ and for all $t>0$. We will say that $K$ is the \emph{kernel} of $T$.
\end{definition}
\medskip
It is clear that a solution admitting the analytic Markov property, will also admit the Markov property. Equation \eqref{markoveq} can be informally understood as $\mu_{t+dt} = (1 + K[\mu_t] d W_t) \mu_t$. So, the function $K$ can be regarded as a transition rule whose randomness comes from the increment $dW_t$. 

In light of equation \eqref{contloc}, the stopping time $T_\mu$ constructed in Section 2 admits the analytic Markov property. In that case, for every $\mu \in \mathcal{M}$, the function $K[\mu]$ is the unique linear function satisfying $ \int_{\RR} K[\mu](x) \mu(dx) = 0$ and $ \int_{\RR} x K[\mu](x) \mu(dx) = 1$. \\

\begin{remark}
The above definitions may be natural in a financial context: when one chooses a market strategy that maximizes the expectation of a certain
quantity, in many cases the optimal strategy need not take the past into account, since the market, $W_t$, is presumably a Markov process.
\end{remark}

This definition gives rise to an entire family of Skorokhod embeddings, as demonstrated by the following proposition. We will omit its proof,
as it follows the same lines as the ones described in Section 2.
\begin{proposition} \label{propmarkov}
Suppose that for a function $K:\mathcal{M} \to \mathcal{MF}$, equation (\ref{markoveq}) has a unique solution for any given $\mu_0=\mu \in \mathcal{M}'$. Then the stopping time induced by equations (\ref{Tmut}) and (\ref{markoveq}) is a Skorokhod embedding if and only if the following conditions hold: \\
(i) $ \int_{\RR} K[\mu](x) \mu(dx) = 0$  for every $\mu \in \mathcal{M}$, \\
(ii) $\int_{\RR} x K[\mu](x) \mu(dx) = 1$  for every $\mu \in \mathcal{M}$, and \\
(iii) $\EE[ \sup \{t|~ Var[\mu_t] > 0\} ] < \infty$. \\ 

\end{proposition}

\begin{remark}
Condition (i) ensures that $\mu_t$ is almost surely a probability measure for all $0<t<T$. Following the same lines as in Lemma \ref{lemdat}, one can see that condition (ii) implies that $\int_{\RR} x \mu_t(dx) = W_t$. Condition (iii) amounts to the fact that $\EE[T] < \infty$.
\end{remark}

\subsection{Bass's embedding}
We now sketch the argument that the embedding introduced by Bass (\cite{Bass}) fits in our framework of embeddings with an analytic Markov property. 

Let $\Phi:\RR \to [0,1]$ be the standard Gaussian cumulative distribution function. For a measure $\nu$, let $\Psi_\nu:\RR \to [0,1]$ be its cumulative distribution function, and let $g_\nu(x) := \Phi^{-1}(\Psi_\nu(x))$. 

Fix a measure $\mu \in \mathcal{M}'$ and a standard Brownian motion $\{B_\tau \}_{\tau=0}^1$. The idea of Bass's embedding is the following: it is clear by definition that $g_\mu^{-1}(B_1) \sim \mu$, therefore, if we define
$$
M_{\tau} = \EE \left .\left [ g_\mu^{-1}(B_1) \right | B_\tau  \right], ~0 \leq \tau \leq 1
$$
then $M_\tau$ is a martingale satisfying $M_1 \sim \mu$. The idea is now to construct a monotone change of time $t(\tau)$ such that $\{M_{\tau(t)} \}_{t=0}^{t(1)}$ has the law of a Brownian motion stopped at $t(1)$, where $\tau(t)$ is the inverse of $t(\tau)$. 

According to a calculation made in Bass's work \cite[Lemma 1]{Bass}, we have that
$$
M_\tau = \int_0^\tau a(\tau, B_s) d B_s
$$
where
$$
a(\tau, y) = \int_{\RR} q_{1-\tau}(z-y) g_\mu^{-1}(z) dz,
$$
with
$$
q_s(y) = -(2 \pi s)^{-1/2} (y/s) e^{-y^2/2s}.
$$
The same calculation also gives that
$$
\frac{d \tau}{dt} = a(\tau(t), B_\tau(t))^{-2}.
$$

Now, suppose that $W_t$ is a standard Brownian motion, into which we wish to embed the measure $\mu$. Then according to Bass's work, we can construct a coupling between $\{B_\tau\}$ and $\{W_t\}$ such that $W_t = M_{\tau(t)}$ for all $t \leq t(1)$. Since $M_1 \sim \mu$, we will have that $W_{t(1)} \sim \mu$, thus $T_\mu = t(1)$ is a solution to Skorokhod problem (the integrability follows from \cite[Lemma 1, Lemma 2]{Bass}).

Now define $\mu_t$ as in \eqref{defmutt}. Since the distribution of $B_1$ conditioned on $B_\tau$ is $N(B_\tau, \sqrt{1-\tau})$, and since $W_{T_\mu} = g_\mu^{-1} \left (B(1) \right )$, we deduce that almost surely for all $0 \leq t \leq T_\mu$, the push forward of $\mu_{\tau(t)}$ under $g_\mu$ has the law $N(B_\tau(t), \sqrt{1-\tau(t)})$. In other words, one has
\begin{equation} \label{pullback}
\frac{d g_\mu^\star \mu_{\tau(t)}}{dx} = \frac{1}{\sqrt{2 \pi (1-t)}} \exp \left (- \frac{ |x - W_t|^2 }{2 (1-t)} \right ),
\end{equation}
where $g^\star \nu$ denotes the push forward of the measure $\nu$ under the map $g$. Using It\^{o}'s formula (see \cite[Lemma 7]{E2} for a detailed calculation), we have
\begin{equation} \label{eqbass}
d \left ( \frac{g_\mu^\star \mu_{\tau(t)} (dx) }{dx}  \right ) = (1-t)^{-1} (x - W_t) \left ( \frac{g_\mu^\star \mu_{\tau(t)} (dx)}{dx}  \right ) d W_t
\end{equation}
for all $x \in \RR$. By using the change of variables $x = g_\mu(y)$ and dividing by the density of $\mu$ with respect to the Lebesgue measure, we get
$$
d \left ( \frac{\mu_{\tau(t)} (dy)}{\mu(dy) } \right ) = (1-t)^{-1} (g_{\mu}(y) - W_t) \left ( \frac{\mu_{\tau(t)} (dy)}{\mu (dy)} \right ) d W_t, ~~ \forall y \in Supp(\mu)
$$
and by integrating with respect to time, we learn that for all measurable $A \subset \RR$, one has
\begin{equation} \label{eqbass1}
\mu_{\tau(t)}(A) = \mu(A) + \int_{0}^t (1-s)^{-1} \int_A (g_{\mu} (x) - W_s) \mu_{\tau(s)}(dx)  d W_s.
\end{equation}
Finally, according to equation \eqref{pullback} we have that $g_\mu \circ g_{\mu(\tau(t))}^{-1}$ pushes forward $N(0,1)$ to $N(W_t, \sqrt{1-t})$ and by the monotonicity of $g_\mu$ one has that 
$$
g_\mu \circ g_{\mu(\tau(t))}^{-1} (y) = W_t + \sqrt{1-t} y, ~ \forall y \in \RR
$$
and thus
$$
g_{\mu_{\tau(t)}} (x) = (1-t)^{-1/2} (g_{\mu}(x) - W_t), ~~ \forall x \in Supp(\mu_{\tau(t)}).
$$
Using this, equation \eqref{eqbass1} becomes
\begin{equation} \label{eqbassfinal}
\mu_{\tau(t)}(A) = \mu(A) + \int_{0}^t (1-s)^{-1/2} \int_A g_{\mu_{\tau(s)}} (x) \mu_{\tau(s)}(dx)  d W_s.
\end{equation}
In view of equation \eqref{markoveq}, we have for Bass's scheme
$$
K_{\mbox{\tiny{Bass}}}[\mu](x) := C_\mu g_\mu (x)
$$
where $C_\mu$ is the unique constant chosen such that condition (ii) of Proposition \ref{propmarkov} holds.

\begin{remark}
Equation \eqref{eqbass} emphasizes a further connection between our construction and the one of Bass. In view of the similarity between this equation and equations \eqref{contloc}, we see that while our embedding, in an infinitesimal time increment, the update of the measure $\mu_t$ corresponds to its multiplication by a linear function, Bass's embedding corresponds to multiplying the push forward of this measure under $g_{\mu_t}$ by a linear function. Nevertheless, since the function $g_{\mu_t}$ changes over time, it becomes clear that those two constructions are generally not the same up to any transformation of time and space.
\end{remark}

\subsection{The Az\'ema-Yor embedding}

When the measure $\mu$ is finitely supported, it turns out that the embedding introduced by Az\'ema and Yor in \cite{AY} satisfies the Markov property as well. 

For a finitely supported measure $\nu$ let $l(\nu) := \min\{Supp(\nu)\}$ be the left-most atom of $\nu$. In case that the cardinality of $Supp(\nu)$ is at least $2$, define $K[\nu]$ to be the unique density supported on $Supp(\nu)$ satisfying: \\
(i) $K[\nu]$ is constant on $Supp(\nu) \setminus \{l(\nu)\}$. \\
(ii) $\int_{\RR} K[\nu](x) \nu(dx) = 0$. \\
(iii) $\int_{\RR} x K[\nu](x) \nu(dx) = 1$. \\

Remark that as long as $Var[\mu_t] > 0$, the cardinality of $Supp(\mu_t)$ has to be at least $2$, therefore equation \eqref{markoveq} has a solution corresponding to the function $K$. We claim that this solution coincides with the Az\'ema-Yor solution to the Skorokhod problem. To see that, fix a measure $\mu$ and let $T$ be the corresponding stopping time constructed according to this embedding scheme. Let $\mu_t$ be defined as in \eqref{defmutt}. Moreover, define $G_t(x) = \frac{\mu_t(dx)}{\mu(dx)}$. By definition of $\mu_t$, we have
$$
G_t(x) = \frac{\PP(W_T = x ~ | \mathcal{F}_t )}{\mu(\{x\} )}, ~~ \forall x \in Supp(\mu)
$$
and consequently $G_t(x)$ is a martingale for all $x$. Moreover, it is evident by the construction of the embedding that $G_t(x)$ is continuous and therefore there exists a predictable process $F_t(x)$ such that
$$
d G_t(x) = F_t(x) G_t(x) d W_t.
$$
Next, according to the construction of this embedding (the reader may refer to \cite[Section 4.1]{Ob} for an accessible description of this construction in the case of finite measures), the following holds: For any $t > 0$, there exists a number $x_t$ such that 
$$
G_t(x) = 0, ~\forall  x \in Supp(\mu) \cap (-\infty, x_t) \mbox{ and } 
$$
$$
G_t(x) = c_t, \forall x \in Supp(\mu) \cap (x_t, \infty) 
$$
where $c_t$ is some (random) constant. In other words, $G_t$ is constant on all but the left-most atom of $\mu_t$. This implies that for all $t$, the function $F_t$ is constant on $Supp(\mu_t) \setminus \{l(\mu_t)\}$, which corresponds exactly to (i) above. Next, since $\mu_t$ is a probability measure for all, $t$, we have that
$$
\int_\RR F_t(x) \mu_t(dx) d W_t = \int_\RR F_t(x) G_t(x) \mu(dx) d W_t 
$$
$$
= d \int_\RR G_t(x) \mu(dx) = d \mu_t(\RR) = 0
$$
and finally since by definition $\int_\RR x \mu_t(dx) = W_t$, we also have
$$
\int_{\RR} x F_t(x) \mu_t(dx) d W_t = d \int_{\RR} x \mu_t(dx) = d W_t
$$
and therefore $\int_\RR x F_t(x) \mu_t(dx) = 1$. We see that conditions (i)-(iii) above hold if $K[\nu]$ is replaced by $F_t$. But, since those conditions uniquely determine the function $F_t$ on $Supp(\mu_t)$, this means that
$$
F_t = K[\mu_t]
$$
which proves that the Az\'ema-Yor embedding must coincide with the construction corresponding to the kernel $K$.

\begin{remark}
When the measure $\mu$ is not finitely supported, the Az\'ema-Yor embedding should still admit the Markov property. However, in general it is not clear to us whether or not it admits the analytic Markov property.
\end{remark}

\subsection{Other known solutions}
It is natural to ask what other known solutions fit under this framework. Let us mention a few examples and questions (several of which were pointed out to us by an anonymous referee, to whom we are grateful). 
\begin{itemize}
\item
The original solution given by Skorokhod \cite{S}, as well is the solution of Hall \cite{Hall} contain extra randomness (meaning that $T$ does not deterministically depend on the path of $\{W_t\}$), and therefore do not admit the analytic Markov property.
\item
The solution of Root \cite{Root} seems to admit the Markov property, as suggested by the work of Loynes (\cite{Loynes}): since the Barrier functions are essentially unique (as shown by Loynes), it follows that the Barrier corresponding to $\mu_t$ should be identical to the barrier corresponding to $\mu$, translated according to $W_t$. However, it is not clear to us whether or not this solution admits the analytic Markov property.
\item
In the solution of Vallois \cite{Vallois}, the fact that $\mu_t$ captures the local times through its support suggests that it may be Markovian. We cannot determine whether or not it is analytically Markovian.
\end{itemize}

\subsection{Another possible Markovian embedding scheme}
It may be natural to consider the following construction: Denote the median of a measure $\nu$ by $med(\nu)$. Consider the equation
$$
F_0 = 1, d F_t = \left ( \mathbf{1}_{ \{x > med(\mu_t) \} } - \mathbf{1}_{ \{x < med(\mu_t) \} } \right )  B_t^{-1} d W_t,
$$
$$
\frac{\mu_t(dx)}{\mu(dx)} = F_t(x) , B_t = \int_{ \{ x > med(\mu_t)\} } x \mu_t(dx) - \int_{ \{x < med(\mu_t)\} } x \mu_t(dx)
$$
note that if the solution of the above equations exists for some initial measure $\mu_0=\mu$, then conditions (i) and (ii) in Proposition \ref{propmarkov} hold if the measure $\mu$ has no atoms.

\begin{question}
Does the above construction induce a Skorokhod embedding for measures with no atoms?
\end{question}

\section{Appendix}
In this appendix we prove that the equation (\ref{stochastic1}) has a unique solution whenever the second moment of $\mu$ is finite. \\ 

We begin with observing that it is enough to prove that almost surely, there exists some $t_0 > 0$ such that the equation has a solution
in the interval $[0,t_0]$. Indeed, $A_\mu (b,c)$ and $a_\mu(b,c)$ are smooth functions any set in which $b$ is bounded away from zero. Since $A_t$
is continuous with respect to $t$ and since $A_0 > 0$, there will necessarily exist some $b' > 0$ such that $b_t \geq b'$ for all $t \geq t_0$ (in other words, the only "problematic" point is $t=0$, since for any $t>0$ the function $\mu_t$, surely has finite exponential moments). \\ 

We argue that there exists a function $c_\mu(\cdot, \cdot)$ satisfying,
\begin{equation} \label{defcmu}
a_\mu(b, c_\mu(a,b)) = a.
\end{equation}
for all $b > 0$ and $a \in Conv(Supp(\mu))$ (the convex hull of the support of $\mu$). Indeed, if we denote that partial derivatives of $a_\mu$ by $a_1(\cdot, \cdot)$ and $a_2(\cdot, \cdot)$ a straightforward calculation gives
\begin{equation} \label{eqa2}
a_2(b,c) = V_{\mu}(b,c)^{-1} \int_{\RR} (x - a_\mu(b,c))^2 e^{c x - \frac{1}{2} b x^2 }  \mu(dx) = A_\mu(b,c) > 0.
\end{equation}
for all $b > 0$. Moreover, it is easy to check that for all $b>0$,
$$
\lim_{c \to -\infty} a_\mu(b,c) = \min Supp(\mu) \mbox{ and } \lim_{c \to \infty} a_\mu(b,c) = \max Supp(\mu).
$$
Observe also that $a_\mu$ is continuous on the domain $b>0$. It follows from inverse function theorem that the function $c_\mu$ exists and is unique. Moreover $c_\mu$ is continuously differentiable in the second argument for all $b>0$ and $a \in Int(Supp(\mu))$. Fix a realization of $W_t$ and define a function
$$
F(b, t) = A_\mu (b, c_\mu(W_t, b))^{-2}
$$
and $F(0,0) = A^{-2}(0,0) = Var[\mu]^{-1}$. We claim that there exists a solution to the equation
\begin{equation} \label{ODE}
\frac{d}{dt} b_t = F(b_t, t), ~~ b_0 = 0
\end{equation}
in an interval $t \in [0,t_0]$. Indeed, since the function $c_\mu(w,b)$ is continuously differentiable with respect to $b$ for $W_t \in Supp(\mu)$ (as indicated above), we have that the function $F(b,t)$ is continuously differentiable with respect to $b$ for all  $(b,w) \in \Omega$ where
$$
\Omega = (0, \infty) \times (0, ~ \sup\{t; W_t \in Supp(\mu) \}).
$$
Consequently, $F(b,t)$ is locally Lipschitz-continuous in the in $b$ on every compact sub-domain of $\Omega$. Moreover, by the continuity of $W_t$ we conclude also that for all $b>0$, the function $F(b,t)$ is continuous with respect to $t$. An application of the Picard-Lindel\"{o}f theorem (see e.g., \cite[Chapter 2]{CL}) establishes the uniqueness and existence of the solution of \eqref{ODE}. \\

Next define $c_t = c_\mu(W_t, b_t)$. Our main goal is to show that $b_t,c_t$ satisfy the equation (\ref{stochastic1}). To that end, we use It\^{o}'s formula to calculate
$$
d c_t =  c_1(W_t, b_t) d W_t + \frac{1}{2} c_{11} (W_t, b_t) dt + c_2(W_t, b_t) \frac{d}{dt} b_t dt.
$$
where $c_{11}(\cdot, \cdot)$ is the second derivative of $c_\mu$ with respect to its first variable. According to the inverse function theorem and using equation (\ref{eqa2}),
$$
c_1(a,b) = A_\mu(b,c_\mu(a,b))^{-1}.
$$
Next, we have
$$
a_1(b,c) = V^{-1}(b, c_\mu(a,b)) \left ( \int_{\RR} \frac{x^2}{2} e^{c_\mu(a,b) x -\frac{1}{2} b x^2} dx a -  \int_{\RR} \frac{x^3}{2} e^{c_\mu(a,b) x -\frac{1}{2} b x^2} dx \right ).
$$
Differentiate equation (\ref{defcmu}) to get
$$
a_1 (b, c_\mu(a,b)) + a_2(b, c_\mu(a,b)) c_2(a,b) = 0
$$
So,
$$
c_2(a,b) = a_1 (b, c_\mu(a,b)) a_2 (b, c_\mu(a,b))^{-1} = 
$$
$$
A_\mu(b,c_\mu(a,b))^{-1} V^{-1}(b, c_\mu(a,b)) \left ( \int_{\RR} \frac{x^2}{2} e^{c_\mu(a,b) x -\frac{1}{2} b x^2} dx a -  \int_{\RR} \frac{x^3}{2} e^{c_\mu(a,b) x -\frac{1}{2} b x^2} dx \right ).
$$
Lastly, we need the second derivative of $c$ with respect to the first variable. One has
$$
c_{11}(a,b) = - a_{22} (b, c_\mu(a,b)) / a_2(b, c_\mu(a,b))^3 = 
$$
$$
A_\mu(b,c_\mu(a,b))^{-3} \left ( V_\mu(b, c\mu(a,b))^{-1} \int_{\RR} (x^3 - a x^2) e^{c_\mu(a,b) x -\frac{1}{2} b x^2} dx + 2 A a  \right ).
$$
We finally get
$$
d c_t =  c_1(W_t, b_t) d W_t + \frac{1}{2} c_{11} (W_t, b_t) dt + c_2(W_t, b_t) \frac{d}{dt} b_t dt = 
$$
$$
A_\mu(b_t, c_\mu(W_t,b_t))^{-1} d W_t + 
$$
$$
\frac{1}{2} A_\mu(b_t,c_t)^{-3} \left ( V_\mu(b_t, c_t)^{-1} \int_{\RR} (x^3 - W_t x^2) e^{c_\mu(W_t, b_t) x -\frac{1}{2} b_t x^2} dx - 2 A_\mu(b_t, c_t) W_t \right ) dt + 
$$
$$
V^{-1}(b_t, c_t) \left ( \int_{\RR} \frac{x^2}{2} e^{c_t x -\frac{1}{2} b_t x^2} dx W_t - \int_{\RR} \frac{x^3}{2} e^{c_\mu(W_t,b_t) x -\frac{1}{2} b_t x^2} dx \right ) A_\mu(b_t, c_\mu(W_t,b_t))^{-3} dt =
$$
$$
A_\mu(b_t, c_t)^{-1} d W_t + W_t A_\mu(b_t, c_t)^{-2} dt.
$$
In view of (\ref{dat}), we conclude that $b_t$ and $c_t$ satisfy (\ref{stochastic1}) in some interval $[0, t_0)$, which proves
the existence and uniqueness of the solution.

\bigskip

{\small \noindent \it e-mail address: roneneldan@gmail.com}

\vfill \hfill \today

\end{document}